\documentclass[11pt]{amsart}

\newtheorem{theorem}{Theorem}[section]
\newtheorem{proposition}[theorem]{Proposition}
\newtheorem{lemma}[theorem]{Lemma}
\newtheorem{claim}[theorem]{Claim}
\newtheorem{fact}[theorem]{Fact}
\newtheorem{cor}[theorem]{Corollary}
 
\newtheorem{definition}[theorem]{Definition}

\newtheorem{conjecture}[theorem]{Conjecture}

\theoremstyle{plain}
\numberwithin{equation}{theorem}

\theoremstyle{remark}
\newtheorem{remark}[theorem]{Remark}
\newtheorem{example}[theorem]{Example}

\newcommand{\F}{{\mathbb F}}

\newcommand{\Kbar}{\overline{K}}

\DeclareMathOperator{\Gal}{Gal}

\DeclareMathOperator{\sep}{sep}
\DeclareMathOperator{\tor}{tor}

\DeclareMathOperator{\CC}{\mathbb{C}_{\infty}}

\DeclareMathOperator{\ord}{ord}

\DeclareMathOperator{\N}{N}

\newcommand{\bP}{{\mathbb P}}
\newcommand{\bZ}{{\mathbb Z}}
\newcommand{\bG}{{\mathbb G}}

\newcommand{\bC}{{\mathbb C}}

\newcommand{\bF}{{\mathbb F}}
\newcommand{\Fq}{\bF_q}
\newcommand{\lra}{\longrightarrow}

\newcommand{\cU}{\mathcal{U}}
\newcommand{\cL}{\mathcal{L}}
\newcommand{\cT}{\mathcal{T}}
\newcommand{\cS}{\mathfrak{S}}

\newcommand{\bK}{\overline{K}}
\newcommand{\hhat}{{\widehat h}}

\newcommand{\Drin}{{\bf \phi}}

\title{Equidistribution and integral points for Drinfeld modules}
\author{D.~Ghioca and T.~J.~Tucker}

\address{
Dragos Ghioca \\
Department of Mathematics\\
McMaster University \\
1280 Main Street West \\ 
Hamilton, Ontario \\
Canada  L8S 4K1 \\
}

\email{dghioca@math.mcmaster.ca}

\address{
Thomas Tucker\\
Department of Mathematics\\
Hylan Building\\
University of Rochester\\
Rochester, NY 14627
}

\email{ttucker@math.rochester.edu}

\begin{document}

\begin{abstract}
  We prove that the local height of a point on a Drinfeld module can
  be computed by averaging the logarithm of the distance to that point over
  the torsion points of the module.  This gives rise to a Drinfeld
  module analog of a weak version of Siegel's integral points
  theorem over number fields and to an analog of a theorem of
  Schinzel's regarding the order of a point modulo certain primes.
\end{abstract}

\maketitle

\section{Introduction}\label{intro} 
In 1929, Siegel (\cite{Siegel}) proved that if $C$ is an affine curve
defined over a number field $K$ and $C$ has at least three points at
infinity, then there are finitely many $K$-rational points with
integral coordinates.  The proof of this famous theorem uses
diophantine approximation along with the fact that certain groups of
rational points are finitely generated; when $C$ has genus greater
than 0, the group in question is the Mordell-Weil group of the
Jacobian of $C$, while when $C$ has genus 0, the group in question is
the group of $S$-units in a finite extension of $K$.

When $C$ is a curve of genus 0, it is simple to give a notion of
integrality that is more flexible than the notion of integral
coordinates.  Let $S$ be a finite set of places of $K$. Viewing
$\bP^1(K)$ as $K$ plus a point at infinity, we say that $\beta \in K$
is $S$-integral with respect to $\alpha \in K$ if at each finite place
$v$ of $K$ outside of $S$, we have $|\alpha - \beta|_v \geq 1$ and
$\min(|\alpha|_v, |\beta|_v) \le 1$ (see Section $2.6$ for a more precise definition of $S$-integrality).  Similarly, we say that $\beta$
is $S$-integral with respect to the point at infinity if $|\beta|_v
\leq 1$ for all points outside of $S$.  Using these conventions, along
with the convention of the usual Weil height $h(\alpha)$ of an
algebraic number $\alpha$, Siegel's method of proof can be described
very easily.  The set of points $\beta$ that are integral with respect
to any two points in $\bP^1_K$ can be realized as a subset of the
group $\Gamma$ of $S$-units in some extension $L$ of $K$.  Since
$\Gamma / m \Gamma$ is finitely generated, if $\alpha$ is another
point (making a total of three) for which there are infinitely many
$\beta_i \in \Gamma$ that are integral with respect to $\alpha$, then
there is some $\gamma \in \Gamma$ such that infinitely many $\gamma +
\beta_i$ have $m$-th roots $\beta_i'$ in $L$.  The fact that the
$\beta_i$ are $S$-integral with respect to $\alpha$ means that
$\sum_{v \in S} - \log |\beta_i - \alpha|_v = h(\beta) + O(1)$ for all
$\beta$.  Then there is an $m$-th root $\alpha'$ of $\alpha + \gamma$,
and infinitely many $m$-th roots $\beta_i' \in L$ of $\beta_i+\gamma$
such that
\begin{equation*}
\begin{split}
\sum_{v \in S} - \log |\beta'_i - \alpha'|_v = \sum_{v \in S} - \log
|\beta_i - \alpha|_v  + O(1) & = h(\beta_i) + O(1) \\
& = m h(\beta'_i) + O(1), 
\end{split}
\end{equation*}
which
violates a theorem from diophantine approximation, such as Roth's
theorem (\cite{Roth}) or Siegel's own result, also found in
\cite{Siegel}.

Now, suppose that the base field of our curves is a function field $K$
over a finite field.  When $C$ is a non-isotrivial curve of positive genus, Siegel's theorem is
still true (see \cite{Voloch-2} and \cite{Voloch-1}), but Siegel's
theorem is not true for the projective line with three points removed.
For example, let $K$ be any function field over a finite field (of
characteristic $p$) and let $x$ be any point in $K$ that is not in the
constant field.  Let $S$ be the set of places for which either $|x|_v
\not = 1$ or $|x-1|_v \not= 1$.  Then any power of $x$ is $S$-integral
with respect to both $0$ and $\infty$ and, since $|x^{p^m} - 1|_v =
|(x-1)^{p^m}|_v = |x-1|_v^{p^m}$ for any positive integer $m$, we see
that $x^{p^m}$ is $S$-integral with respect to $1$ as well.  Note that
the diophantine approximation results over function fields are much
weaker than Roth's (or Siegel's or Dyson's) theorem (see \cite{VolSur}
or \cite{Osgood}, for example).

In this paper we will prove an analog of Siegel's theorem for cyclic
$\phi$-submodules of $\bG_a(K)$ (under the action of a Drinfeld module
$\phi$). In some sense this is surprising in light of the
counterexample above for $\bZ$-submodules of $\bG_m(K)$.  On the other
hand, there are many other famous theorems over number fields that
have analogs in characteristic $p$ in the context of Drinfeld modules.
For example, Ghioca \cite{IMRN} and Scanlon \cite{Scanlon} proved
Mordell-Lang and Manin-Mumford statements, respectively, for Drinfeld
modules; these had both been conjectured by Denis (\cite{Denis-conjectures}).
Bosser formulated and proved (\cite{Bosser-thesis}) a variant of the
Bogomolov conjecture for Drinfeld modules (unpublished, but the main
ingredient of his result is contained in \cite{Bosser-compositio}).
Ghioca (\cite{Mat.Ann}) proved a Drinfeld module analog of Szpiro,
Ullmo, and Zhang's equidistribution theorem (\cite{SUZ}) for abelian
varieties.  In this paper we explore further aspects of
equidistribution in the context of Drinfeld modules and demonstrate a
connection between equidistribution and integrality of points.

We will prove the following theorem.

\begin{theorem}\label{integral:1}
 Let $K$ be a finite extension of $\F_q(t)$. Let $\phi:\mathbb{F}_q[t]\rightarrow K\{\tau\}$ be a Drinfeld module and let $S$ be a finite set of places of $K$.  If $\beta\in K$ is not
  torsion for $\phi$, then there are finitely many $Q\in\Fq[t]$ such that $\Drin_Q(\beta)$
  is $S$-integral for the point $0$.
\end{theorem}

We will derive Theorem \ref{integral:1} as a consequence of Theorem
\ref{easiest}, which states that the local height of a point $\beta$
can be computed by averaging the function $\log |x - \beta|_v$ over
the periodic points of $\Drin$.  Theorem \ref{easiest} can be regarded
as an equidistribution theorem, since local canonical heights can
typically be computed by integrating $\log |x - \beta|_v$ against an
invariant measure (see Section \ref{further} for more details).
Theorem \ref{easiest} also gives rise to Theorem \ref{precise
  order}, which is a result about the order of the reduction of a
point modulo various primes. Determining the order of the reduction of
a nontorsion point at various primes was also studied in the context
of abelian varieties (see \cite{Pink-abelian}).

Let us now give a brief outline of the paper.  We begin with some
preliminaries about canonical heights on Drinfeld modules in Section
\ref{notation}.  Following that we prove Theorem \ref{easiest}. The main ingredient for Theorem~\ref{easiest} is Corollary~\ref{alternative local height computation} which provides a new way of defining local canonical heights for a Drinfeld module. For the places lying over the place at infinity for the Drinfeld module (see Section~\ref{notation} for a more precise definition of these places), this is proved using linear
forms in logarithms (as in \cite{Bosser}); for the ``finite''
places, this is proved using an elementary method described in
Lemmas~\ref{small-x} and \ref{L:nice trick}.  Theorem \ref{integral:1}
follows easily from Theorem \ref{easiest}.  The proof of Theorem
\ref{precise order} requires a somewhat lengthy counting argument,
which is detailed in Section \ref{schinzel}.  In Section
\ref{further}, we discuss further generalizations, including analogs
of work of Silverman (\cite{SilSiegel}), Baker/Ih/Rumely (\cite{BIR}),
and Szpiro/Tucker (\cite{ST}).  The most general question we pose may
be viewed as a possible Siegel's integral points theorem for arbitrary
finitely generated $\phi$-submodules of $\bG_a(\bK)$. In a recent paper \cite{our_siegel} we answer positively our question under one mild technical hypothesis and also assuming a natural conjecture for linear forms in logarithms at finite places for a Drinfeld module (which is the analog of Bosser's result \cite{Bosser} for the linear forms in logarithms at infinite places).

\noindent{\it Acknowledgments.}  We would like to thank M. Baker and
V. Bosser for many interesting conversations.

\section{Notation}
\label{notation}

\subsection{Drinfeld modules}
We begin by defining a Drinfeld module.  Let $p$ be a prime and let
$q$ be a power of $p$. Let $A:=\mathbb{F}_q[t]$, let $K$ be a finite extension of $\mathbb{F}_q(t)$, and let $\Kbar$ be an
algebraic closure of $K$. We let $\tau$ be the Frobenius on
$\mathbb{F}_q$, and we extend its action on $\Kbar$.  Let $K\{\tau\}$
be the ring of polynomials in $\tau$ with coefficients from $K$ (the
addition is the usual addition, while the multiplication is the
composition of functions).

A Drinfeld module is a morphism $\Drin:A\rightarrow K\{\tau\}$ for
which the coefficient of $\tau^0$ in $\Drin(a)=:\Drin_a$ is $a$ for
every $a\in A$, and there exists $a\in A$ such that $\Drin_a\ne
a\tau^0$. The definition given here represents what Goss \cite{Goss}
calls a Drinfeld module of ``generic characteristic''.

We note that usually, in the definition of a Drinfeld module, $A$ is the ring of functions defined on a projective nonsingular curve $C$, regular away from a closed point $\eta\in C$. For our definition of a Drinfeld module, $C=\mathbb{P}^1_{\mathbb{F}_q}$ and $\eta$ is the usual point at infinity on $\mathbb{P}^1$. On the other hand, every ring of regular functions $A$ as above contains $\mathbb{F}_q[t]$ as a subring, where $t$ is a nonconstant function in $A$. 

A point $\alpha$ is \emph{torsion} for the Drinfeld module action if
and only if there exists $Q\in A\setminus\mathbb{F}_q$ such that
$\Drin_Q(\alpha)=0$. The monic polynomial $Q$ of minimal degree which
satisfies $\phi_Q(\alpha)=0$ is called the \emph{order} of $\alpha$.
Since each polynomial $\phi_Q$ is separable, the torsion submodule
$\phi_{\tor}$ lies in the separable closure $K^{\sep}$ of $K$.
 
For every field extension $K\subset L$, the Drinfeld module $\Drin$
induces an action on $\mathbb{G}_a(L)$ by $a*x:=\Drin_a(x)$, for each
$a\in A$. We call \emph{$\phi$-submodules} subgroups of
$\mathbb{G}_a(\overline{K})$ which are invariant under the action of
$\phi$. As shown in \cite{Poonen}, $\mathbb{G}_a(K)$
is a direct sum of a finite torsion $\phi$-submodule with a free
$\phi$-submodule of rank $\aleph_0$.

We will often use sums of the form $\sum_{P|Q}$ over divisors $P\in A$
of a fixed polynomial $Q\in A$. These sums will always be taken over
the \emph{monic} divisors $P$ of $Q$.

\subsection{Valuations and Weil heights}
Let $M_{\mathbb{F}_q(t)}$ be the set of places on $\Fq(t)$.
We denote by $v_{\infty}$ the place in $M_{\Fq(t)}$ such that $v_{\infty}(\frac{f}{g})=\deg(g)-\deg(f)$ for
every nonzero $f,g\in A=\Fq[t]$. We let $M_K$ be the set of valuations on $K$. Then $M_K$ is a set of valuations which satisfies a
product formula (see \cite[Chapter 2]{Serre-Mordell_Weil}). Thus
\begin{itemize}
\item for each nonzero $x\in K$, there are finitely many $v\in M_K$
such that $|x|_v\ne 1$; and \\
\item for each nonzero $x\in K$, we have $\prod_{v\in M_K} |x|_v=1$.
\end{itemize}
We may use these valuations to define a
Weil height for each $x \in K$ as
\begin{equation}\label{weil}
h(x) = \sum_{v \in M_K} \log \max (|x|_v,1).
\end{equation}

We may also extend the Weil height to all of $\bK$ in a coherent
way (see Chapter $2$ of \cite{Serre-Mordell_Weil} for a detailed discussion of the
construction of heights on $\bK$).  More specifically, if $L$ is an
extension of $K$, then we can define a set $M_L$ of absolute values
$| \cdot |_{w}$ on $L$ such that
\begin{equation}\label{defectless}
|y|_v = \prod_{\substack{w|v\\w\in M_{L}}}| y |_{w}
\end{equation}
Note that we write $w | v$ when $| \cdot |_w$ restricts to some
power of $| \cdot |_v$ on $K$, but that $| \cdot |_w$ is {\it not}
in general equal to $| \cdot |_v$ on $K$ (see \cite[pages
9-11]{Serre-Mordell_Weil}). We use these $| \cdot |_w$ to define a
height for any $y$ in a finite extension $L$ of $K$ by
\begin{equation}\label{h(y)}
h(y) = \sum_{w \in M_L} \log \max(| y |_w, 1).
\end{equation}
It follows from \eqref{defectless}, that for any $x \in K$ and any
extension $L$ of $K$, we have
\begin{equation}\label{gen-weil}
\sum_{v \in M_K} \log \max(|x|_v, 1) = \sum_{w \in M_L} \log
\max(| x |_w, 1).
\end{equation}
Thus the height on all of $\bK$ restricts to the original Weil height
\eqref{weil} on $K$.  More generally, the definition of $h(y)$ in
\eqref{h(y)} does not depend on our choice of the field $L$ containing
$y$.  For more details, we refer the reader to \cite[Chapter
2]{Serre-Mordell_Weil}.

\begin{definition}
  Let $L$ be a finite extension of $K$. Each place in $M_L$ that lies
  over $v_{\infty}$ is called an infinite place. Each place in $M_L$
  that does not lie over $v_{\infty}$ is called a finite place.
\end{definition}

\subsection{Canonical heights}
Let $\Drin:A\rightarrow K\{\tau\}$ be a Drinfeld module of \emph{rank}
$d$ (i.e. the degree of $\phi_t$ as a polynomial in $\tau$ equals
$d$). The canonical height of $\beta\in\bK$ relative to $\Drin$ (see \cite{Denis}) is
defined as
$$\hhat(\beta) = \lim_{n \to \infty}
\frac{h(\Drin_{t^n}(\beta))}{q^{nd}}.$$ Denis \cite{Denis} showed that $\alpha$ is a torsion point for $\Drin$
if and only if $\hhat(\alpha)=0$.

For every finite extension $L$ of $K$, and for every $w\in M_L$, we let the local canonical height of $\beta\in L$ at $w$ be
\begin{equation}\label{poon-def}
\hhat_w(\beta) = \lim_{n \to \infty} \frac{\log
  \max(|\Drin_{t^n}(\beta)|_w, 1)}{q^{nd}}.
\end{equation}
It is clear that $\sum_{w\in M_L} \hhat_w(\beta) = \hhat(\beta)$.

We will use the notion of \emph{Galois conjugate} for points
$\beta\in\bK$. Even though this notion is classical, for the sake of
completeness, we define it here. If $\sigma:\bK\rightarrow\bK$ is a
field automorphism fixing $K$, then $\beta^{\sigma}:=\sigma(\beta)$ is
called a Galois conjugate of $\beta$. 

\subsection{Completions and filled Julia sets}

By abuse of notation, we let $\infty\in M_K$ denote any place extending the
place $v_{\infty}$.  We let $K_{\infty}$ be
the completion of $K$ with respect to $|\cdot|_{\infty}$. We let
$\overline{K_{\infty}}$ be an algebraic closure of $K_{\infty}$. We
let $\CC$ be the completion of $\overline{K_{\infty}}$. Then $\CC$ is
a complete, algebraically closed field. Note that $\CC$ depends on our
choice for $\infty\in M_K$ extending $v_{\infty}$. However, each time
we will work with only one such place $\infty$, and so, there will be
no possibility of confusion.

Next, we define the \emph{$w$-adic filled Julia set} $J_{\phi,w}$
corresponding to the Drinfeld module $\phi$ and to each place $w$ of
$M_L$ (where $L$ is any finite extension of $K$).  Let $\mathbb{C}_w$ be the completion of an algebraic closure
of $L_w$.  Then $| \cdot |_w$ extends to a unique absolute value on
all of $\mathbb{C}_w$.  The set $J_{\phi,w}$ consists of all
$x\in\mathbb{C}_w$ for which $\{|\phi_Q(x)|_w\}_{Q\in A}$ is bounded.
It is immediate to see that $x\in J_{\phi,w}$ if and only if
$\{|\phi_{t^n}(x)|_w\}_{n\ge 1}$ is bounded.

One final note on absolute values: as noted above, each place $v\in M_K$ extends to a unique absolute value $|\cdot |_v$ on all of $\mathbb{C}_v$.  We fix
an embedding of $i: \bK  \lra \mathbb{C}_v$.  For $x \in\bK$, we denote
$| i(x) |_v$ simply as $|x|_v$, by abuse of notation.

\subsection{The coefficients of $\Drin_t$}
Each Drinfeld module is isomorphic to a Drinfeld module for which all
the coefficients of $\Drin_t$ are integral at all the places in $M_K$
which do not lie over $v_{\infty}$. Indeed, we let $B\in
\Fq[t]$ be a product of all (the finitely many)
irreducible polynomials $P\in\Fq[t]$ with the property
that there exists a place $v\in M_K$ which lies over the place $(P)\in M_{\Fq(t)}$, and there exists a coefficient of $\phi_t$
which is not integral at $v$. Let $\gamma$ be a sufficiently large
power of $B$. Then $\psi:A\rightarrow K\{\tau\}$ defined by
$\psi_Q:=\gamma^{-1}\Drin_Q\gamma$ (for each $Q\in A$) is a Drinfeld
module isomorphic to $\Drin$, and all the coefficients of $\psi_t$ are
integral away from the places lying above $v_{\infty}$. Hence, from
now on, we assume that all the coefficients of $\Drin_t$ are integral
away from the places lying over $v_{\infty}$. It follows that for
every $Q\in A$, all coefficients of $\phi_Q$ are integral away from
the places lying over $v_{\infty}$.

\subsection{Integrality and reduction}
\begin{definition}
\label{S-integral definition}
For a finite set of places $S \subset M_K$ and $\alpha\in\bK$, we say
that $\beta \in\bK$ is $S$-integral with respect to $\alpha$ if for
every place $v\notin S$, and for every pair of morphisms $\sigma,\tau:\bK\rightarrow\bK$ the following
are true:
\begin{itemize}
\item if $|\alpha^{\tau}|_v\le 1$, then $|\alpha^{\tau}-\beta^{\sigma}|_v\ge 1$.
\item if $|\alpha^{\tau}|_v >1$, then $|\beta^{\sigma}|_v\le 1$.
\end{itemize}
\end{definition}
For more details about the definition of $S$-integrality, we refer the reader to \cite{BIR}.

We also introduce the notation $\overline{\beta}$ for the reduction at
a place $w\in M_L$ of a point $\beta\in L$, which is integral at $w$.
When $\phi$ has \emph{good reduction} at $w$ (i.e., for each nonzero
$Q\in\Fq[t]$, all coefficients of $\phi_Q$ are $w$-adic integers and
the leading coefficient of $\phi_Q$ is a $w$-adic unit), we denote by
$\overline{\phi}$ the reduction of $\phi$ at $w$. Note that when
$\phi$ has good reduction, $\overline{\phi}$ is a well-defined Drinfeld module on $\bG_a(k_w)$ of the same rank as $\phi$, where $k_w$ is the residue field at $w$. When
$\phi$ does not have good reduction at $w$, we say that it has {\em
  bad reduction} at $w$.  We observe that each place lying above
$v_{\infty}$ is a place of bad reduction for $\phi$.

The following fact is proved in an equivalent form also in
\cite{locleh} (see Lemma $4.13$).
\begin{fact}
\label{good places torsion}
Let $w$ be a place of
good reduction for $\phi$, and let $\alpha\in\bC_w$ be a torsion point for $\phi$. Then $\alpha$ is a $w$-adic integer.
\end{fact}

\begin{proof}
  Assume $\alpha$ is not a $w$-adic integer. Let $Q\in\Fq[t]$ be the
  order of $\alpha$. Assume $\phi_Q=\sum_{j=0}^l b_j\tau^j$. Then for
  every $0\le j<l$,
$$|b_l\alpha^{q^l}|_w=|\alpha^{q^l}|_w>|\alpha^{q^j}|_w\ge |b_j\alpha^{q^j}|_w.$$
Hence $|\phi_Q(\alpha)|_w=|\alpha^{q^l}|_w>1$, which contradicts $\phi_Q(\alpha)=0$.
\end{proof}

\section{Proofs of main theorems}\label{main}

We continue with the notation as in Section~\ref{notation}. The main goal of this section is to prove the following theorem.
Then Theorem~\ref{integral:1} will follow easily as a corollary. 

\begin{theorem}\label{easiest}
  For any nontorsion $\beta \in \bK$ and any nonzero
  irreducible $F \in K[Z]$ such that $F(\beta) = 0$, we have
\begin{equation} \label{ht-comp}
\begin{split}
  (\deg F) \hhat(\beta) =\sum_{v\in M_K} \lim_{\deg Q \to \infty}
  \frac{1}{{q^{d \deg Q}}}  \sum_{\Drin_Q(y) = 0} \log
  | F(y) |_v.
\end{split}
\end{equation}
\end{theorem}

Later (see Theorem~\ref{also easiest}), we will prove a similar result
as in Theorem~\ref{easiest} valid also for torsion points $\beta$. Our
proof of Theorem~\ref{easiest} will go through a series of lemmas and
propositions which ultimately will yield a new way to compute local
heights (see Corollary~\ref{alternative local height computation}).
We will use in our argument a Proposition proved by Bosser (see
Th\'{e}or\`{e}me $1.1$ of \cite{Bosser}) as well as a Lemma due to
Taguchi (see Lemma $(4.2)$ of \cite{Taguchi}).

\begin{proposition}\label{from Bosser}(Bosser)
Let $\infty$ be a place lying over $v_{\infty}$. Let $\exp_{\phi}$ be the exponential map associated to $\phi$ at $\infty$ (see \cite{Goss}).
Let $y_1, \dots, y_n$ be elements of $\CC$ such that
$\exp_{\Drin}(y_i) \in \bG_a(\bK)$.  Then there exists a (negative) constant $C_1$ (depending on $\phi$, $y_1,\dots,y_n$) such
that for any $P_1, \dots, P_n \in \Fq[t]$, either $P_1y_1+\dots+P_ny_n=0$ or 
$$ \log | P_1 y_1 + \dots + P_n y_n|_{\infty} \geq C_1 \max_{1\le i\le n}(\deg(P_i)\log
\deg(P_i)).
  $$
\end{proposition}

\begin{remark}
  In the translation of Th\'{e}or\`{e}me $1.1$ of \cite{Bosser} to our
  Proposition~\ref{from Bosser}, we use the fact that for each
  polynomial $P\in\F_q[t]$, its height is its degree.
\end{remark}

\begin{lemma}\label{from Taguchi}(Taguchi)
  Let $\exp_{\Drin}$ be the associated exponential map to $\Drin$ with respect to a fixed place $\infty$ lying over $v_{\infty}$. Let $\cL$ be the lattice of all elements 
  $u \in\CC$
  such that $\exp_{\Drin}(u) = 0.$ There exists an $A$-basis $\omega_1,
  \dots, \omega_d$ for $\cL$ such that for any polynomials $P_1,
  \dots, P_d \in \Fq[t]$, we have
$$ |P_1 \omega_1 + \dots +P_d \omega_d|_\infty = \max_{1\le i\le d} |P_i \omega_i|_\infty.$$ 
\end{lemma}
\begin{remark}
The elements $\omega_i$ from Lemma~\ref{from Taguchi} form a basis for $\cL$ of ``successive minima'' (as defined in \cite{Taguchi}).
\end{remark}

The following result appears in a similar form in \cite{Goss} (see
Proposition $4.14.2$).
\begin{lemma}\label{small ball}
With the notation as in Lemma~\ref{from Taguchi}, there exists a constant $C_2 > 0$ such that $\exp_{\Drin}$ induces an
isomorphism from the metric space $B_1=\{ z \in \bG_a(\bC_{\infty}) \; \mid \; |z|_{\infty} \leq C_2
\}$ to itself. 
\end{lemma}
\begin{proof}
As shown in \cite{Goss},
$$\exp_{\Drin}(z)=z\prod_{\substack{\lambda\in\cL\\ \lambda\ne
    0}}(1-\frac{z}{\lambda}).$$
Because $\cL$ is discrete in
$|\cdot|_{\infty}$, there exists $C_2>0$ such that for every nonzero
$\lambda\in\cL$, we have $|\lambda|_{\infty}> C_2$. Thus, for each
such $\lambda$, if $|z|_{\infty}\le C_2$, then
$|1-\frac{z}{\lambda}|_{\infty}=1$. Therefore, the ball $$B_1:=\{ z
\in \bG_a(\CC) \; \mid \; |z|_{\infty} \leq C_2\}$$
is mapped by
$\exp_{\Drin}$ into itself (preserving the distance). Moreover,
$\exp_{\Drin}$ is injective on $B_1$ because $\cL\cap B_1=\{0\}$.
Finally, $\exp_{\Drin}:B_1\rightarrow B_1$ is surjective as can be
easily seen from considering the entire function
$f_{\alpha}(z):=\exp_{\Drin}(z)-\alpha$ (for each $\alpha\in B_1$).
Then $f_{\alpha}$ has a zero on $B_1$ (see the first slope of its
associated Newton polygon).
\end{proof}

\begin{proposition}
\label{P:first iteration}
Let $L$ be a finite extension of $K$, and let $\infty\in M_L$ be a place lying above $v_{\infty}$. Let $\beta\in L$ be any nontorsion point for the Drinfeld module $\Drin$ that is in $J_{\phi,\infty}$. Then
$$ \lim_{\deg Q \to \infty} \frac{\log |\Drin_Q(\beta)|_{\infty}}{q^{d \deg Q}} =
0. $$
\end{proposition}
\begin{proof}
  Let $\epsilon$ be any real number greater than zero.  Since $\beta$
  is in the filled Julia set of $\Drin$, we know that
  $|\phi_Q(\beta)|_\infty$ is bounded, so
$$ \frac{\log |\phi_Q(\beta)|_{\infty}}{q^{d \deg Q}} < \epsilon$$
when $\deg Q$ is sufficiently large.

Now, if $\frac{\log |\Drin_Q(\beta)|_{\infty}}{q^{d \deg Q}} < -
\epsilon$ for $\deg(Q)$ sufficiently large, then we have
$|\Drin_Q(\beta)|_{\infty} < C_2$ where $C_2$ is as in Lemma
\ref{small ball}.  We will use Proposition~\ref{from Bosser} to derive
a contradiction.  First, we let $u\in\CC$ such that $\exp_{\phi}(u)=\beta$. We also choose a
basis $\omega_1,\dots,\omega_d$ of ``successive minima'' as in
Lemma~\ref{from Taguchi}. Throughout our argument, $u,\omega_1,\dots,\omega_d$ are fixed.

Because $|\Drin_Q(\beta)|_{\infty}<C_2$, then there exists $y$ such
that $|y|_{\infty}=|\Drin_Q(\beta)|_{\infty}$ and
$\exp_{\Drin}(y)=\Drin_Q(\beta)$. Since $Qu$ satisfies
$\exp_{\Drin}(Qu)=\Drin_Q(\beta)$, there exist $P_1,\dots,P_d\in A$
such that
$$y=Qu+P_1\omega_1+\dots+P_d\omega_d.$$
Moreover, because of
Lemma~\ref{small ball}, there exists no $y'$ with
$$|y'|_{\infty}<|\Drin_Q(\beta)|_{\infty}<C_2$$
such that $\exp_{\Drin}(y')=\Drin_Q(\beta)$. Thus, $|Qu|_{\infty}\ge
|y|_{\infty}$. Therefore, we have
$|P_1\omega_1+\dots+P_d\omega_d|_{\infty}\le |Qu|_{\infty}$ (because
otherwise $|y|_{\infty}>|Qu|_{\infty}$).

Thus, by Lemma \ref{from Taguchi} we see that
 \begin{equation}
 \label{equation-Taguchi}
 \max_i |P_i \omega_i|_\infty \leq |Q|_\infty \cdot |u|_\infty.
 \end{equation}
Taking logarithms in \eqref{equation-Taguchi}, we obtain
$$\max_i \deg(P_i)\le
\deg(Q)+\log|u|_{\infty}-\min_i\log|\omega_i|_{\infty}.$$
Let
$C_3:=\log|u|_{\infty}$.  Hence, $\max_i \deg(P_i) \leq \deg Q +
C_4$ (the constant $C_4$ depends only on $C_3$ and $\phi$, as the
$(\omega_i)_{i=1}^d$ is a fixed basis of successive minima for the
corresponding lattice of $\phi$ at $\infty$). Thus,
Proposition~\ref{from Bosser} gives
\begin{equation}
\label{E:Bosser}
\begin{split}
  \log |\Drin_Q(\beta)|_{\infty} = \log |y|_\infty & = \log|Qu + P_1 \omega_1 + \dots + P_d
  \omega_d|_{\infty} \\
  & \geq C_1 (\deg Q + C_4)\log(\deg Q + C_4).
\end{split}
\end{equation}
Note that we know $y\ne 0$, because $\phi_Q(\beta)\ne 0$ (and hence,
we can apply the inequality from Lemma~\ref{from Bosser}).  For large
$Q$, this means that $\frac{\log |\Drin_Q(\beta)|_\infty}{q^{d \deg
    Q}} \geq - \epsilon$, which completes our proof.
\end{proof}

\begin{cor}
\label{C:first iteration}
For any nontorsion $\beta \in \bK$ and for any place $\infty\in M_{K(\beta)}$ lying above $v_{\infty}$, we have
$$ \lim_{\deg Q \to \infty} \frac{\log |\Drin_Q(\beta)|_{\infty}}{q^{d \deg Q}} =
\hhat_{\infty}(\beta). $$

\end{cor} 
\begin{proof}
  If $\beta$ is in the filled Julia set, then this is
  Proposition~\ref{P:first iteration}. If $\beta$ is not in the filled
  Julia set, then for large $n$, we have
  $$|\Drin_{t^n}(\beta)|_\infty >
  \max(|\Drin_{t^m}(\beta)|_\infty,1)$$
  for all $m < n$. So, for $Q$ of
  large degree we have
$$|\Drin_{Q}(\beta)|_\infty = |\Drin_{t^{\deg Q}}(\beta)|_\infty =
\max (| \Drin_{t^{\deg Q}}(\beta)|_\infty,1).$$
  Our claim then follows from
  the definition of the local canonical height with respect to
  $\Drin$.
\end{proof}

Let $L$ be a finite extension of $K$ containing $\beta$. Now we deal with the places $w$ of $L$ which do not lie over $v_{\infty}$. We
recall that for such $w$, all coefficients of $\phi_Q$ are integral at
$w$ for each $Q\in A$.

If $\beta\notin J_{\phi,w}$, then $|\phi_Q(\beta)|_w$ is unbounded
as $\deg Q \to \infty$.  Then, by the definition of the local
canonical height (see also the proof of Corollary~\ref{C:first
  iteration}), we have
\begin{equation}
\label{E:new height 1}
\hhat_w(\beta)=\lim_{\deg(Q)\rightarrow\infty}\frac{\log|\Drin_Q(\beta)|_w}{q^{d\deg(Q)}}.
\end{equation}
We will show that when $\beta\in J_{\phi,w}$, then the limit above is
equal to $0$.
\begin{proposition}
\label{P:first finite place iteration}
Let $L$ be a finite extension of $K$, and let $w\in M_L$ be a place not lying over $v_\infty$. Let $\beta\in L$ be a nontorsion point which is also in $J_{\phi,w}$. Then
$$\lim_{\deg(Q)\rightarrow\infty}\frac{\log|\Drin_Q(\beta)|_w}{q^{d\deg(Q)}}=0.$$
\end{proposition}

\begin{proof}
  Because $\beta$ is in the $w$-adic filled Julia set, then
  $|\phi_Q(\beta)|_w$ is bounded. Hence, for every $\epsilon>0$, there
  exists an integer $N^{+}_{\epsilon}$ such that if
  $\deg(Q)>N^{+}_{\epsilon}$, then
  $$\frac{\log|\phi_Q(\beta)|_w}{q^{d\deg(Q)}}<\epsilon.$$
  It remains
  to show that for every $\epsilon>0$ there exists an integer
  $N^{-}_{\epsilon}$ such that if $\deg(Q)>N^{-}_{\epsilon}$, then
\begin{equation}
\label{N-}
\frac{\log|\phi_Q(\beta)|_w}{q^{d\deg(Q)}}>-\epsilon.
\end{equation}
For proving \eqref{N-}, we start with a general lemma about valuations.
\begin{lemma}
\label{small-x}
Let $F(X) = b_n X^n + \dots + b_1 X$ be a polynomial such that $|b_i|_w
\leq 1$ for all $i$.  Then $|F(x)|_w = |b_1 x|_w$ whenever $|x|_w <
|b_1|_w$.  In particular, $|F(x)|_w \leq |x|_w$.
\end{lemma}
\begin{proof}[Proof of Lemma~\ref{small-x}.]
For any $i \geq 2$, we have 
$$|b_i x^i|_w \leq |b_i x^2|_w < |b_i b_1 x|_w \leq |b_1 x|_w.$$
Applying the ultrametric inequality thus yields
$$ |F(x)|_w = |b_n x^n + \dots + b_1 x|_w = |b_1 x|_w .$$
\end{proof}

The following result is an immediate corollary of the fact that all coefficients of $\phi_Q$ (for $Q\in A$) are $w$-adic integers.
\begin{cor}
\label{C:ideal}
   For $0 <\epsilon \le 1$, the set of polynomials $Q \in A$ such that
  $|\phi_{Q}(\beta)|_w < \epsilon$ is an ideal in $A$.
\end{cor}

We continue the proof of Proposition~\ref{P:first finite place iteration}. The place $w$ restricts to a place on $\Fq(t)$ which corresponds to an irreducible polynomial $P\in A$.

Hence, for each polynomial $Q\in A$, we have
$|Q|_w=|P|_w^{\ord_P(Q)}$, where $P^{\ord_P(Q)}$ is the largest
power of $P$ dividing $Q$. The following lemma will finish the proof
of Proposition~\ref{P:first finite place iteration}.

\begin{lemma}
\label{L:nice trick}
  There exists a positive constant $C_w$
  such that
  $$
  |\phi_Q(\beta)|_w \geq C_w |P|_w^{\ord_P(Q)}$$
for all $Q \in A$.  
\end{lemma}
\begin{proof}
  If there is no $F \in A$ such that $|\phi_{F}(\beta)|_w < |P|_w$,
  then we are done.  Otherwise, let $G\in A$ generate the ideal of all
  polynomials $F$ such that $|\phi_{F}(\beta)|_w < |P|_w<1$, and let
  $C_w := |\Drin_G(\beta)|_w$ (note that $\phi_G(\beta)\ne 0$ because
  $\beta$ is nontorsion).
  
  Let $Q$ be any element in the ideal generated by $G$.  We may write
  $Q = FG$ for a polynomial $F\in\Fq[t]$. We expand $F$ out
  $P$-adically, i.e. we write
$$ F = R_n P^n + \dots + R_1 P + R_0 $$
where each $R_i$ is a polynomial of degree less than $\deg P$.  Since
$|R_i|_w = 1$ for each $i$, we have
$$
|(\phi_{R_i}\phi_G) (\beta)|_w = C_w < |P|_w $$
by Lemma \ref{small-x}.
Moreover, also by Lemma~\ref{small-x}, we have 
$$
\left| (\phi_P \phi_{R_i} \phi_G) (\beta)
\right|_w = C_w |P|_w,$$
and more generally, by induction, we have 
$$
\left| (\phi_{P^i}\phi_{R_i}\phi_G) (\beta)
\right|_w = C_w |P|^i_w,$$
for any $i$.  Thus, letting $f$ be the smallest $i$ such that $R_i
\not= 0$, we have 
$$ |\phi_Q(\beta)|_w = |(\phi_{R_fP^f}\phi_G)(\beta)|_w = C_w
|P|_w^f.$$
Since $f \leq \ord_P(Q)$, this completes our proof of Lemma~\ref{L:nice trick}.  
\end{proof}
Using Lemma~\ref{L:nice trick}, we obtain 
$$\log|\phi_Q(\beta)|_w\ge \log(C_w)+\deg(Q)\log|P|_w,$$
which proves \eqref{N-} and thus concludes the proof of
Proposition~\ref{P:first finite place iteration}.
\end{proof}

The following corollary is an immedate consequence of
Proposition~\ref{P:first finite place iteration} (see also
Corollary~\ref{C:first iteration}).
\begin{cor}
\label{alternative local height computation}
Let $L$ be a finite extension of $K$ and let $\beta\in L$ be a nontorsion point. Then for every $w\in M_{L}$, $$\hhat_w(\beta)=\lim_{\deg(Q)\rightarrow\infty}\frac{\log|\phi_Q(\beta)|_w}{\deg(\phi_Q)}.$$
\end{cor}

We are now ready to prove Theorem~\ref{easiest}.  

\begin{proof}[Proof of Theorem~\ref{easiest}.]
  Using the product formula applied to the leading coefficient of $F$,
  we see that it suffices to prove our result under the assumption
  that $F$ is monic.  Write $F = \prod_{i=1}^n (Z - \theta_i)$.  Then
\begin{equation}
\label{radacini}
\begin{split}
  & \sum_{v \in M_K} \lim_{\deg Q\rightarrow\infty} \frac{1}{q^{d\deg
      Q}} \sum_{\phi_Q(y)=0} \log|F(y)|_v \\ 
  & = \sum_{v \in M_K} \lim_{\deg Q\rightarrow\infty}
  \frac{1}{q^{d\deg Q}} \sum_{\phi_Q(y)=0} \sum_{i=1}^n\log|y-\theta_i|_v \\
  & = \sum_{v\in M_K} \lim_{\deg Q\rightarrow\infty}
  \frac{1}{q^{d\deg(Q)}} \sum_{i=1}^n\log|\phi_Q(\theta_i)|_v\\
  & \quad \quad - \sum_{v\in M_K} \lim_{\deg Q\rightarrow\infty}
  \frac{n \log|\gamma_Q|_{v}}{q^{d \deg Q}} ,
\end{split}
\end{equation} 
where $\gamma_Q$ is the leading coefficient of $\phi_Q$.  Using
induction, it is easy to see that $\log|\gamma_Q|_v = {\frac{q^{d\deg
      Q}-1}{q^d-1}} \log |c a_d|_v$, where $a_d$ is the leading
coefficient of $\phi_t$ and $c$ is the leading coefficient of $Q$
(which is in $\Fq$).  Passing to the limit, we find then that
$$\lim_{\deg Q\to\infty}\frac{\log |\gamma_Q|_v}{q^{d\deg Q}}=\frac{\log |a_d|_v}{q^d-1}.$$
Applying the product formula to $a_d$ we thus see that the sum
$$\sum_{v\in M_K} \lim_{\deg Q\rightarrow\infty} \frac{\log|\gamma_Q|_{v}}{q^{d \deg
    Q}} $$
vanishes.  Hence, it will suffice to show that
\begin{equation}\label{nearly}
  (\deg F) (\hhat(\beta)) =  \sum_{v \in M_K}\lim_{\deg Q \to \infty}  \sum_{i=1}^n
  \frac{\log  |\phi_Q(\theta_i)|_v}{q^{d \deg Q}}.
\end{equation}

Now, let $L:= K(\beta)$, and let $p^t$ (for some $t\ge 0$) be the
inseparable degree of $L/K$. Then $n=p^ts$, where $s$ is the separable
degree of $L/K$. Moreover, the $\theta_i$ above (for $1\le i\le n$)
split into $s$ multisets, each of them consisting of $p^t$ identical
elements (while different multisets have different elements).
Therefore we may choose distinct $\theta_1,\dots,\theta_s$ among the
$\theta_i$ above. Let $K_1:=K^{1/p^t}$, and let $L_1$ be the
compositum (inside $\bK$) of $K_1$ and $L$. Then $L_1/K_1$ is a
separable extension of degree $s$. Moreover, the irreducible
polynomial $F_1$ of $\beta\in L_1$ over $K_1$ equals $F^{1/p^t}$.
Dividing both sides of \eqref{nearly} by $p^t$, we obtain

\begin{equation}\label{nearly-separable}
   (\deg F_1) (\hhat(\beta)) = \sum_{v \in M_{K_1}} \lim_{\deg Q \to \infty}  \sum_{i=1}^s
  \frac{\log |\phi_Q(\theta_i)|_v}{q^{d \deg Q}}.
\end{equation}

Note that for \eqref{nearly-separable} we used the fact that above each place $v\in M_K$ lies a unique place of $K_1$ (which we also denote by $v$). Therefore we are left to prove \eqref{nearly} in the case $\beta$ generates a separable extension. Thus from now on we assume $L/K$ is a separable extension of degree $n$.

Following \cite[page
11]{Serre-Mordell_Weil}, each absolute value of some $x\in L_w$ (for $w|v$) takes the
form $|\N_{L_w/K_v}(x)|_v^{1/[L:K]}$.  Writing $\N_{L_w/K_v}(x)$ as
$\iota_{w,1} (x) \cdots \iota_{w,m_w}(x)$, where $\iota_{w,1},
\dots, \iota_{w, m_w}$ are the embeddings of $L_w$ into $\bC_v$ (and $m_w=[L_w:K_v]$), we
have 
$$ | x |_w^{[L:K]} = \prod_{j=1}^{m_w} |\iota_{w,j}(x)|_v.$$
For each $\theta_i$ in $\bC_v$, there is exactly one pair $(w,j)$ such
that $\iota_{w,j}(\beta) = \theta_i$. So, using that $\phi_Q(\theta_i)=\iota_{w,j}(\phi_Q(\beta))$ (because $\phi_Q$ has coefficients in $K$), we sum over all $L_w$ and obtain
\begin{equation*}
\begin{split}
  \sum_{i=1}^n \log |\phi_Q(\theta_i)|_v 
  & = [L:K] \sum_{w \mid v} \log | \phi_Q(\beta) |_w \\
 \end{split}
\end{equation*}  
Passing to the limit, and using the fact
that $[L:K] = (\deg F)$, we thus obtain
\begin{equation*}
\begin{split}
  \sum_{v \in M_K} \lim_{\deg Q \to \infty} \sum_{i=1}^n \frac{\log
    |\phi_Q(\theta_i)|_v}{q^{d \deg Q}} & = \sum_{w
    \in M_L}\lim_{\deg Q\to\infty} (\deg F)\cdot \frac{\log  |\phi_Q(\beta) |_w}{{q^{d \deg Q}}}   \\
  & =(\deg F)\sum_{w\in M_L}\hhat_w(\beta)\\
  & = (\deg F) \hhat(\beta)
\end{split}
\end{equation*}
as desired.
\end{proof}

We can extend the conclusion of Theorem~\ref{easiest} and prove the following result.
\begin{theorem}\label{also easiest}
  For any $\beta \in \bK$ and any nonzero
  irreducible $F \in K[Z]$ such that $F(\beta) = 0$, we have
\begin{equation} \label{also ht-comp}
\begin{split}
  (\deg F) \hhat(\beta) =\sum_{v\in M_K} \lim_{\deg Q \to \infty}
  \frac{1}{{q^{d \deg Q}}} \sum_{\substack{\Drin_Q(y)
      = 0\\ y\ne \beta}} \log | F(y) |_v.
\end{split}
\end{equation}
\end{theorem}

\begin{proof}
  The difference from Theorem~\ref{easiest} is that in
  Theorem~\ref{also easiest} we deal also with torsion points $\beta$.
  This is why we exclude from the inner summation of \eqref{also
    ht-comp} the case $y=\beta$. Theorem~\ref{also easiest} for
  $\beta\in\phi_{\tor}$ follows exactly as Theorem~\ref{easiest} once
  we prove the following Proposition.
\begin{proposition}
\label{fact-torsion}
Assume $\beta\in\bK$ is a torsion point. Then 
$$\sum_{v\in M_K}\lim_{\deg Q\rightarrow\infty}\frac{1}{q^{d\deg
    Q}} \sum_{\substack{\phi_Q(y)=0\\ 
    y\ne\beta}}\log|y-\beta|_{v}=0.$$
\end{proposition}

\begin{proof}[Proof of Proposition~\ref{fact-torsion}.]
Let $v\in M_K$. For each $Q\in A$ we let $\gamma_Q$ be the leading
coefficient of $\phi_Q$. Also, let $a_d$ be the leading coefficient of
$\phi_t$. We will prove that
\begin{equation}
\label{primul termen}
\lim_{\deg Q\rightarrow\infty} \frac{1}{q^{d\deg Q}}
 \sum_{\substack{\phi_Q(y)=0\\
    y\ne\beta}} \log|y-\beta|_{v}= -\frac{\log|a_d|_v}{q^d-1}.
\end{equation}
Using \eqref{primul termen} and the product formula for $a_d$ we will
then conclude the proof of Proposition~\ref{fact-torsion}.

We proceed to proving \eqref{primul termen}. Let $Q\in A$ be a non-constant polynomial. There are two cases:

\emph{Case 1.} $\phi_Q(\beta)=0$

Then $$\prod_{\substack{\phi_Q(y)=0\\ 
    y\ne\beta}}(\beta-y)=\frac{\left(\phi_Q\right)'(\beta)}{\gamma_Q}.$$
But the derivative $\left(\phi_Q\right)'$ is identically equal to $Q$.
Since $|\gamma_Q|_v=|a_d|_v^{\frac{q^{d\deg Q}-1}{q^d-1}}$, we have
\begin{equation*}
\begin{split}
  \lim_{\deg Q\rightarrow\infty} \frac{1}{q^{d\deg
      Q}} \sum_{\substack{\phi_Q(y)=0\\
      y\ne\beta}}\log|y-\beta|_{v} & =\lim_{\deg
    Q\rightarrow\infty}\frac{\log|Q|_v-\frac{q^{d\deg
        Q}-1}{q^d-1}\log|a_d|_v}{q^{d\deg Q}}\\ &
  =-\frac{\log|a_d|_v}{q^d-1},
\end{split}
\end{equation*}
as desired.

\emph{Case 2.} $\phi_Q(\beta)\ne 0$

Then 
\begin{equation}
\label{gamma-q}
\prod_{\substack{\phi_Q(y)=0\\ y\ne\beta}}(\beta-y)=\prod_{\phi_Q(y)=0}(\beta-y)=\frac{\phi_Q(\beta)}{\gamma_Q}.
\end{equation}
However, $\phi_Q(\beta)$ is one of the finitely many nonzero torsion points of $\phi$ of order dividing the order of $\beta$. Therefore, taking logarithms of \eqref{gamma-q} and using $|\gamma_Q|_v=|a_d|_v^{\frac{q^{d\deg Q}-1}{q^d-1}}$, we conclude the proof of \eqref{primul termen}.
\end{proof}

Arguing as in the proof of Theorem~\ref{easiest} we deduce the
remaining case of Theorem~\ref{also easiest} (i.e.
$\beta\in\phi_{\tor}$) from Proposition~\ref{fact-torsion}.
\end{proof}

We are now ready to prove Theorem \ref{integral:1}.
\begin{proof}[Proof of Theorem \ref{integral:1}.]
  Since $\beta$ is in $K$, the polynomial $F$ satisfied by $\beta$ in
  the statement of Theorem~\ref{easiest} may be taken to be $F(Z) = Z
  - \beta$.  Let
\begin{equation}\label{Drin-poly}
\Drin_t(X) = a_dX^{q^d} + a_{d-1}X^{q^{d-1}}+ \dots +a_1X^q+ t X.
\end{equation}
Then the leading coefficient $\gamma_Q$ of $\Drin_Q$ equals
$ca_d^{\frac{q^{d\deg(Q)}-1}{q^d-1}}$, where $c\in\Fq$ is the leading
coefficient of $Q$. We observe that $|\gamma_Q|_v\ne 1$ if and only if
$|a_d|_v\ne 1$. Then, after taking products followed by logarithms, we see
that
\begin{equation}\label{trick}
\sum_{\Drin_Q(y) = 0} \log |y-\beta|_v = \log |\Drin_Q(\beta)|_v-\log|\gamma_Q|_v.
\end{equation}
Now, suppose that there are infinitely many $Q$ such that
$\Drin_Q(\beta)$ is $S$-integral for $0$.  Let $\cT$ be the set of all
places in $S$ along with all the places $v$ for which $|\beta|_v > 1$
and all the places of bad reduction for $\phi$.  For any $Q$, at any
$v$ outside of $\cT$, we have $|y|_v\le 1$ when $\Drin_Q(y) = 0$ (see
Fact~\ref{good places torsion}).  Since $|\beta|_v\le 1$, then $\log
|y-\beta|_v \leq 0$ for $v\notin\cT$.

When $\Drin_Q(\beta)$ is $S$-integral for 0, it is certainly
$\cT$-integral as well. So, for all $v$ outside $\cT$, we have
$\sum_y\log |y-\beta|_v \geq 0$, by \eqref{trick} (note that
$|\gamma_Q|_v=1$ and $|\phi_Q(\beta)|_v\ge 1$ because $\phi_Q(\beta)$ is $\cT$-integral).  Hence, for these $Q$, we have $\log |y-\beta|_v
= 0$ for all $v$ outside $\cT$, and for all $y$ such that
$\phi_Q(y)=0$.  Since $\cT$ is finite, we can interchange the limit
and the sum in \eqref{ht-comp} and apply the product formula, which
yields $\hhat(\beta) = 0$. This is a contradiction, since $\beta$ is
nontorsion.
\end{proof}

\begin{remark}
  It is worth noting that our alternative way of computing local
  heights (see Corollary~\ref{alternative local height computation})
  also holds for $\beta\in\mathbb{C}_w$ if $w$ does not lie over
  $v_{\infty}$. For example, our key Lemma~\ref{L:nice trick} holds
  for every nontorsion point in $\mathbb{C}_w$, not only for algebraic
  points (i.e. for points in $\bK$). However,
  Corollary~\ref{alternative local height computation} cannot be
  extended to all $\beta\in\CC$ as shown by the following example.
  This is the case because Bosser's result (our Lemma~\ref{from
    Bosser}) does not hold for linear forms in logarithms of
  transcendental points in $\CC$ (i.e points not in $\bK$).
\end{remark}

\begin{example}
We construct an example of a point $\beta\in\CC$ for which the limit $\lim_{\deg Q\to\infty}\frac{\log|\phi_Q(\beta)|_{\infty}}{q^{d\deg Q}}$ does not exist.

We work with the simplest Drinfeld module (but a similar construction
works for any Drinfeld module as the reader will easily see). So, we
let $\phi :\Fq[t]\rightarrow\Fq(t)\{\tau\}$ be the Carlitz module, i.e. $\phi_t(x) = tx + x^q$. By abuse of notation, we denote by $\infty$ the place $v_{\infty}$ of $\Fq(t)$.  The Carlitz module has
associated to it a $1$-dimensional lattice at $\infty$ which is
spanned by an element called $\zeta$, say.

Let $(d(n))_{n\ge 1}$ be a sequence of positive integers which is very
rapidly increasing, i.e. $d(n+1)$ is much larger than $q^{d(n)}$. For
the sake of concreteness, we will use the sequence $d(n)$ defined
recursively, by letting $d(1):=q$ and for every $n\ge 1$,
$d(n+1):=q^{q^{d(n)}}$.

Let $\alpha := \zeta\cdot \sum_{n\ge 1} \frac{1}{t^{d(n)}} \in \CC$.
Clearly $\alpha$ is a (very rapidly) convergent power series.

Let $\beta := \exp_{\phi}(\alpha)$, where $\exp_{\phi}$ is the usual
exponential map of $\phi$ at $\infty$. This is the transcendental
point (i.e. not in $\bK$) which will give us our counterexample because
$|\phi_{t^{d(n)}}(\beta)|_{\infty}$ will be much too rapidly
decreasing.

Indeed, $\phi_{t^{d(n)}}(\beta) = \exp_{\phi}(t^{d(n)}\alpha)$ and
using the fact that $\exp_{\phi}$ kills everything in
$\Fq[t]\cdot\zeta$, we conclude

$$\phi_{t^{d(n)}}(\beta) = \exp_{\phi}\left(\zeta\cdot \sum_{m>n}
  \frac{1}{t^{d(m)-d(n)}}\right).$$

Now, because our sequence $(d(m))_{m\ge 1}$ was chosen to be rapidly
increasing, then $|\frac{1}{t^{d(n+1)-d(n)}}|_{\infty}$ is
sufficiently small (for sufficiently large $n$) and so, $\zeta\cdot
\sum_{m>n} \frac{1}{t^{d(m)-d(n)}}$ is in a sufficiently small ball
around $0$ on which $\exp_{\phi}$ preserves the metric (see
Lemma~\ref{small ball}). Moreover, 
$$\left|\zeta \cdot \sum_{m>n} \frac{1}{t^{d(m)-d(n)}}\right|_{\infty} =
|\zeta|_{\infty}\cdot \left|\frac{1}{t^{d(n+1)-d(n)}}\right|_{\infty}.$$
Therefore, we have
\begin{equation}
|\phi_{t^{d(n)}}(\beta)|_{\infty} = |\zeta|_{\infty} \cdot \left|\frac{1}{t^{d(n+1)-d(n)}}\right|_{\infty}
= |\zeta|_{\infty} \cdot e^{-(d(n+1)-d(n))},
\end{equation}
so $\log |\phi_{t^{d(n)}}(\beta)|_{\infty} = \log|\zeta|_{\infty} -
(d(n+1)-d(n))$. Since $d(n+1)-d(n)$ is much larger than $q^{d(n)}$
(because $d(n+1)=q^{q^{d(n)}}$), we see that $\log
|\phi_{t^{d(n)}}(\beta)|_{\infty}$ is much smaller than $ - q^{d(n)}$.
More precisely,
$$\lim_{n\rightarrow\infty}\frac{d(n+1)-d(n)}{q^{d(n)}}=\lim_{n\rightarrow\infty}\frac{d(n+1)}{q^{d(n)}}=\lim_{n\rightarrow\infty}\frac{q^{q^{d(n)}}}{q^{d(n)}}=+\infty,$$
showing that we do not have a finite limit for
$$\frac{\log |\phi_Q(\beta)|_{\infty}}{q^{\deg(Q)}}\text{ as $\deg Q\to\infty$.}$$
\end{example}

\section{An analog of a theorem of Schinzel}\label{schinzel}

Using our findings we are able to prove the following result, which is
similar to a result of Schinzel \cite{Schinzel} for primitive divisors
of $B^n-C^n$ in number fields. It came to our attention that
independently, Hsia \cite{Hsia} proved a similar Schinzel statement
for Drinfeld modules. We thank him for pointing out the argument for
our Lemma~\ref{more precise}.

 \begin{theorem}
\label{precise order}
Let $K$ be a finite extension of $\Fq(t)$. Let $\beta\in K$ be a
nontorsion point for $\phi$, and let $S$ be a finite set of places in
$M_K$. Then there exists a positive integer $N$ such that for all
$Q\in\Fq[t]$ of degree at least $N$, there exists $v\in M_K\setminus
S$ such that
\begin{enumerate}
\item $|\phi_Q(\beta)|_v<1$; and
\item for every nonzero polynomial $P\in\Fq[t]$ that divides $Q$ and
  has smaller degree than $Q$, we have $|\phi_P(\beta)|_v\ge 1$.
\end{enumerate}
\end{theorem}
 
We begin by proving a precise result combining Lemmas~\ref{small-x}
and \ref{L:nice trick}.

\begin{lemma}
\label{more precise}
Let $v\in M_K$ be a place lying above the irreducible polynomial
$P\in\Fq[t]$. Let $l:=\deg P$. For any $x\in K$, if 
$$|x|_v < |P|_v^{\frac{1}{q^l-1}}<1,$$
then for every $Q\in\Fq[t]$, we
have $|\phi_Q(x)|_v=|Q x|_v\le |x|_v$.
\end{lemma}

\begin{proof}
  If $Q\in\Fq[t]$ is coprime with $P$, then $|Q|_v=1$. Thus
  $|x|_v<|Q|_v$ and so, using Lemma~\ref{small-x}, we conclude that
  $|\phi_Q(x)|_v=|Qx|_v$.
  
  Let $\phi_P=\sum_{i=0}^{dl}b_i\tau^i$. We know that $b_0=P$ (by the definition of $\phi$). If
\begin{equation}
\label{bound}
|x|_v<\min_{1\le i\le dl}\left( \left| \frac{P}{b_i} \right|_v\right)^{\frac{1}{q^i-1}},
\end{equation}
then for every $1\le i\le dl$, we have $|Px|_v>|b_ix^{q^i}|_v$; so,
$|\phi_P(x)|_v=|Px|_v$. We will show next that $|x|_v <
|P|_v^{\frac{1}{q^l-1}}$ implies \eqref{bound}. First we prove that
$|b_i|_v\le |P|_v$ for every $1\le i\le l-1$.

If $l=1$, then the above claim is vacuously true. Therefore, assume
$l>1$. We will use induction on $i$ to show that $|b_i|_v\le |P|_v$
for every $1\le i\le l-1$. Since $\phi_P\phi_t = \phi_t\phi_P$,
writing $\phi_t=\sum_{i=0}^d a_i\tau^i,$ and noting that the coefficient
of $\tau$ is the same in $\phi_t\phi_P$ as it is in $\phi_P\phi_t$, we
obtain
\begin{equation}
\label{1-coefficient}
a_0b_1+a_1b_0^q=b_0a_1+b_1a_0^q.
\end{equation}
Since $b_0=P$ and $a_0=t$, \eqref{1-coefficient} gives rise to the
equality
$$|b_1|_v\cdot |t^q-t|_v=|a_1|_v\cdot |P^q-P|_v\le 1\cdot
|P|_v=|P|_v.$$
Since we assumed $l>1$, we have $|t^q-t|_v=1$. Hence,
the equation above yields $|b_1|_v\le |P|_v$, as desired.  Proceeding
by induction, we assume now that $|b_j|_v\le |P|_v$ for all $1\le
j<i$. Equating coefficients of $\tau^i$ in both $\phi_t\phi_P$ and
$\phi_P\phi_t$, we obtain
\begin{equation*}
\label{i-coefficient}
a_0 b_i + \sum_{j>1} a_j b_{i-j}^{q^j}=b_ia_0^{q^i}+\sum_{j<i}b_ja_{i-j}^{q^j}.
\end{equation*}
Since $|b_j|_v\le |P|_v$ for all $1\le j<i$ by assumption (and
$|a_j|_v\le 1$ for every $j$), we see that $|b_i|_v\cdot
|t^{q^i}-t|_v\le |P|_v.$ But $t^{q^i}-t$ is divisible only by
irreducible polynomials in $\Fq[t]$ of degree at most $i$ (since
$t^{q^i - 1} - 1$ splits completely over a degree $i$ extension of
$\Fq$).  Therefore, $|t^{q^i}-t|_v=1$ (because $i<l$), so $|b_i|_v\le
|P|_v$, as desired.

Since $|b_i|_v\le |P|_v<1$ for every $1\le i\le l-1$, and $|b_i|_v\le
1$ for every $i\ge l$, we have
\begin{equation*}
|P|_v^{\frac{1}{q^l-1}}\le \min_{1\le i\le
  dl}\left( \left| \frac{P}{b_i} \right|_v\right)^{\frac{1}{q^i-1}}.
\end{equation*}
So, if $|x|_v\le |P|_v^{\frac{1}{q^l-1}}$, then \eqref{bound} holds and
we have$|\phi_P(x)|_v=|Px|_v$.  Moreover, $|Px|_v=|P|_v\cdot
|x|_v<|P|_v$ and so, using Lemma~\ref{small-x}, we see that
$$|\phi_{P^2}(x)|_v = |\phi_P(\phi_P(x))|_v = |P\phi_P(x)|_v =
|P^2x|_v.$$
Similarly, an easy induction shows that for every $i\ge 0$, we have
$|\phi_{P^i}(x)|_v = |P^ix|_v$. For any (nonzero) $Q\in\Fq[t]$, there
exists $i\ge 0$ and $R\in\Fq[t]$ coprime with $P$ such that $Q=P^iR$.
Using $|\phi_R(x)|_v=|x|_v$, we thus obtain
$$|\phi_Q(x)|_v=|\phi_{P^i}(\phi_R(x))|_v = |P^i|_v\cdot
|\phi_R(x)|_v=|Q|_v\cdot |x|_v,$$
and our proof is complete.
\end{proof}

\begin{cor}
\label{util}
Let $x\in K$. For all but finitely many places $v\in M_K$, if
$\pi_v\in K$ is an uniformizer at $v$, and $|x|_v\le |\pi_v|_v$, then
$|\phi_Q(x)|_v=|Qx|_v$ for every $Q\in\Fq[t]$. More precisely, if
$|x|_v<1$ and $v$ lies over a place of $\Fq(t)$ corresponding to an
irreducible polynomial of degree larger than
$\log_q\left([K:\Fq(t)]+1\right)$, then $|\phi_Q(x)|_v=|Qx|_v$ for
every $Q\in\Fq[t]$.  
\end{cor}

\begin{proof}
  Let $v\in M_K$ lie over the irreducible polynomial $P\in\Fq[t]$, and
  let $l=\deg P$. By Lemma~\ref{more precise}, if $|x|_v <
  |P|_v^{\frac{1}{q^l-1}}$, then $|\phi_Q(x)|_v=|Qx|_v$ for every
  $Q\in\Fq[t]$. But $|P|_v\ge |\pi_v|_v^{[K:\Fq(t)]}$.  Therefore, if
  $q^l-1>[K:\Fq(t)]$, then
  $$|P|_v^{\frac{1}{q^l-1}} > |\pi_v|_v.$$
  Thus, if $|x|_v\le
  |\pi_v|_v$, then (by Lemma~\ref{more precise}) $|\phi_Q(x)|_v=|Qx|_v$ for every $Q\in\Fq[t]$. Finally, note that there are finitely
  many polynomials in $\Fq[t]$ having bounded degree.
\end{proof}

We will also need the following technical result about local heights
for our proof of Theorem~\ref{precise order}.

\begin{lemma}
\label{positive local heights}
Let $\beta\in K$ be a nontorsion point for the Drinfeld module $\phi$.
Then there exist positive constants $C_0$ and $N_0$ (depending on
$\phi$ and $\beta$) such that for every place $v\in M_K$ for which
$\hhat_v(\beta)>0$, and for every $Q\in A=\Fq[t]$ of degree larger
than $N_0$, we have
$$q^{d\deg(Q)}\hhat_v(\beta)- C_0 \le \log|\phi_Q(\beta)|_v \le
q^{d\deg(Q)}\hhat_v(\beta)+C_0.$$
\end{lemma}

\begin{proof}
  Because we allow the constants $C_0$ and $N_0$ from the conclusion of
  Lemma~\ref{positive local heights} to depend on $\beta$, it suffices
  to prove our lemma for each (of the finitely many places) $v\in M_K$
  for which $\hhat_v(\beta)>0$.

Let $\phi_t=\sum_{i=0}^d a_i\tau^i$. For a fixed place $v\in M_K$ for
which $\hhat_v(\beta)>0$, we define
$$M_v:=\max \left\{ \left(\frac{|a_i|_v}{|a_d|_v}\right)^{\frac{1}{q^d-q^i}}\text{ : }
0\le i<d\right\}\cup\left\{\frac{1}{|a_d|_v^{\frac{1}{q^d-1}}} \right\}.$$
Let $N_0$ be
a positive integer such that $|\phi_{t^{N_0}}(\beta)|_v>M_v$ (we can
find such $N_0$ because $\left(|\phi_{t^n}(\beta)|_v\right)_n$ is
unbounded). Let $\gamma:=|\phi_{t^{N_0}}(\beta)|_v$.  The definition
of $M_v$ yields
$$|\phi_{t^{N_0+1}}(\beta)|_v=|\phi_{t}(\phi_{t^{N_0}}(\beta))|_v=|a_d|_v\cdot\gamma^{q^d}>\gamma>M_v.$$
A simple induction shows that for every $n\ge 1$, we have
$$|\phi_{t^{N_0+n}}(\beta)|_v=|a_d|_v^{\frac{q^{dn}-1}{q^d-1}}\gamma^{q^{nd}}.$$
Hence, 
$$\hhat_v(\beta)=\lim_{n\rightarrow\infty}
\frac{\log|\phi_{t^{N_0+n}}(\beta)|_v}{q^{d(N_0+n)}}=\frac{\log|a_d|_v}{q^{dN_0}(q^d-1)}+\frac{\log\gamma}{q^{dN_0}},$$
which means that
$$\log|\phi_{t^{N_0+n}}(\beta)|_v=q^{d(N_0+n)}\hhat_v(\beta)-\frac{\log|a_d|_v}{q^d-1}.$$
Moreover, for every $Q\in A$ of degree $N_0+n$ (for $n\ge 1$), we have
$$\log |\phi_Q(\beta)|_v = \log
|\phi_{t^{N_0+n}}(\beta)|_v=q^{d(N_0+n)}
\hhat_v(\beta)-\frac{\log|a_d|_v}{q^d-1},$$
which concludes the proof
of Lemma~\ref{positive local heights}.
\end{proof}

The following lemma is the key to proving Theorem~\ref{precise order}.
\begin{lemma}\label{one v}
With the hypothesis from Theorem~\ref{precise order}, let $Q \in \Fq[t]$ be a non-constant monic polynomial.  
Let $\cT$ be the set of places $v\in M_K$ that satisfy the following properties
\begin{enumerate}
\item if $P\in\Fq[t]$ is an irreducible polynomial for which $|P|_v<1$, then $\deg P>\deg_q\left([K:\Fq(t)]+1\right)$;
\item $\beta$ is a $v$-adic integer;
\item $v$ is a place of good reduction for $\phi$;
\item either $|\phi_Q(\beta)|_v=1$, or there is some $P \not= Q$ in $\Fq[t]$ such that $P | Q$ and
  $|\Drin_P(\beta)|_v < 1$;
\item $v\notin S$.
\end{enumerate}
Then
\begin{equation}
 \sum_{v \in \cT} \sum_{y \text{ has order }Q} \log
 |\beta - y|_{v} \geq - \deg Q.
\end{equation}
\end{lemma}

\begin{proof}
  Since 
  $$\sum_{\substack{v \in M_K \\v \text{ does not lie over }
      v_{\infty}}} \log |Q|_v = -\log|Q|_{v_{\infty}}= -\deg Q$$
  and no
  place $v\in\cT$ lies over $v_{\infty}$, it will suffice to show that
  for each $v\in\cT$, we have
\begin{equation}\label{one v eq}
 \sum_{\text{$y$ has order $Q$}} \log |\beta - y|_{v} \geq  \log |Q|_v
\end{equation}

If $|\phi_Q(\beta)|_v=1$, then for every $y$ such that $\phi_Q(y)=0$,
we have $|y-\beta|_v=1$. Indeed, Fact~\ref{good places torsion}
implies $|y|_v\le 1$ for all torsion points $y$, because all $v\in\cT$
are places of good reduction for $\phi$. Moreover, $|\beta-y|_v\le 1$
for $v\in\cT$ and $y\in\phi_{\tor}$ due to $(ii)$.  Furthermore, the
leading coefficient $\gamma_Q$ of $\phi_Q$ is a $v$-adic unit. Thus
$$\log|\phi_Q(\beta)|_v=\sum_{\phi_Q(y)=0}\log|y-\beta|_v.$$
Hence, if $|\phi_Q(\beta)|_v=1$, then indeed $|y-\beta|_v=1$ for all $y$ of order dividing $Q$. Therefore \eqref{one v eq} holds in this case because $|Q|_v\le 1$ for $v\in\cT$.

Thus, we may assume that $|\phi_Q(\beta)|_v<1$.  Let $P_0$ be the
smallest degree monic polynomial dividing $Q$ such that
$|\Drin_{P_0}(\beta)|_v < 1$.  Then, since $|y-\beta|_v \leq 1$ for
each torsion point $y$ and since the leading coefficient of
$\phi_{P_0}$ is a $v$-adic unit, we have
\begin{equation}
\label{sleek 1}
\begin{split}
  \sum_{\substack{\phi_P(y) = 0 \\
      \text { for $P | Q$ with $P \not = Q$}}} \log |\beta - y|_{v}
  \le \sum_{\phi_{P_0}(y) = 0} \log |\beta - y|_{v} =
  \log|\phi_{P_0}(\beta)|_{v}.
\end{split}
\end{equation}
Using again that the leading coefficient of $\phi_Q$ is a $v$-adic unit, we obtain
\begin{equation}
\label{sleek 0}
\sum_{\phi_Q(y) = 0} \log |\beta
- y|_{v} = \log |\Drin_Q(\beta)|_v.
\end{equation}
Since $v\in\cT$ satisfies $(i)$ and $|\phi_{P_0}(x)|_v <1 $, we can
use the more precise claim of Corollary~\ref{util} and derive
\begin{equation}
\label{sleek 2}
\log|\Drin_Q(\beta)|_v =\log\left| \frac{Q}{P_0} \right|_v +\log |\Drin_{P_0}(\beta)|_v \geq \log|Q|_v +
\log|\Drin_{P_0}(\beta)|_v.
\end{equation}
Equations \eqref{sleek 1}, \eqref{sleek 0}, and \eqref{sleek 2} yield
\eqref{one v eq}, which finishes the proof of Lemma~\ref{one v}.
\end{proof}

We define the M\"{o}bius function $\mu$ on the multiplicative set of
all monic polynomials in $\mathbb{F}_q[t]$ by
$$\mu(1)=1,$$
$$\mu(Q_1Q_2\dots Q_n)=(-1)^n,$$
if $Q_1,\dots,Q_n$ are distinct irreducible, non-constant polynomials, and
$$\mu(f)=0\text{ if $f$ is not squarefree.}$$
\begin{lemma}
\label{inclusion-exclusion}
With the notation as in Lemma~\ref{one v}, then for each $v\in\cT$, we
have
$$\sum_{y\text{ has order }Q}\log|y-\beta|_{v} =
\sum_{P|Q}\mu\left(\frac{Q}{P}\right)\log|\phi_P(\beta)|_v.$$
\end{lemma}

\begin{proof}
  Since the leading coefficient of $\phi_P$ (for each nonzero
  polynomial $P$) is a $v$-adic unit for each place $v\in\cT$, then
  for each $P|Q$, we have
\begin{equation}
\label{roots}
\sum_{\Drin_P(y)=0}\log|y-\beta|_{v} = \log|\phi_P(\beta)|_v.
\end{equation}
Using \eqref{roots} and the principle of inclusion and exclusion applied to the set of all $y$ such that $\phi_Q(y)=0$ (by counting them with respect to their corresponding orders $P|Q$), we obtain the result of Lemma~\ref{inclusion-exclusion}.
\end{proof}

We are now ready to prove Theorem~\ref{precise order}.
\begin{proof}[Proof of Theorem~\ref{precise order}.]
  After dividing $Q$ by its leading coefficient (which is a constant
  in $\Fq$), we may assume $Q$ is monic.  We will argue by
  contradiction.
    
  We begin by dividing the places of $K$ into three sets.  We let
  $\cS$ denote the set of all places $v$ such that $\hhat_v(\beta) >
  0$, let $\cT$ be as in Lemma \ref{one v}, and let $\cU$ denote the
  remaining places.  Note that $\cS$ and $\cU$ are both finite,
  because $\cT$ contains all but finitely many places in $M_K$ (here
  we are using the assumption that the conclusion of our
  Theorem~\ref{precise order} fails, because then condition $(iv)$ in
  Lemma~\ref{one v} is satisfied by all but finitely many places).
  Note that $S\subset\cS\cup\cU$ because by assumption $(v)$ from
  Lemma~\ref{one v}, the set $\cT\cap S$ is empty.

We begin by dealing with the places in $\cS$.  As proved in
Lemma~\ref{positive local heights}, there exists a positive integer
$N_0$ and a positive constant $C_0$ such that for all $Q\in\Fq[t]$ of
degree larger than $N_0$, and for all $v\in\cS$, we have 
\begin{equation}
\label{E:exact height}
q^{d\deg(Q)} \hhat_v(\beta) - C_0\le \log|\phi_Q(\beta)|_v \le
q^{d\deg(Q)} \hhat_v(\beta) + C_0.
\end{equation}
Therefore, from now on, we will always assume $\deg(Q)> N_0$.   

Next, we treat the places in $\cU$.  As proved in \eqref{E:Bosser} and
in Lemma~\ref{L:nice trick}, there exist positive constants $C_1$,
$C_2$, and $C_3$ depending only on $\phi$, $\beta$ and $\cU$ such that
for all $v\in \cU$, then
\begin{equation}
\label{E:almost height 0}
-C_1\deg(Q)\log\deg(Q)-C_2<\log|\phi_Q(\beta)|_v<C_3.
\end{equation}
The right hand side of \eqref{E:almost height 0} is guaranteed by the
fact that $\hhat_v(\beta)=0$ (because only $\cS$ contains places $v$ for which $\hhat_v(\beta)>0$).  Summing over all the $v$ in $\cU$, we
see that there are constants $C_4$, $C_5$, and $C_6$ such that
\begin{equation}\label{upper lower}
-C_4 \deg(Q) \log\deg(Q) - C_5 < \sum_{v \in \cU} \log |\phi_Q(\beta)|_v < C_6.
\end{equation}
Since there are finitely many polynomials $P$ of degree at most equal to $N_0$, we see from equations \eqref{E:exact height} and \eqref{upper lower} that there are constants $C_7$ and $C_8$ such that
\begin{equation} 
\begin{split}
  q^{d\deg P} \hhat(\beta) - C_4 \deg(Q) \log\deg(Q)-C_7 & \le
  \sum_{v \in \cS \cup \cU} \log |\phi_P(\beta)|_v  \\ & \le q^{d \deg
    P} \hhat(\beta) + C_8,
\end{split}
\end{equation}
for all $P|Q$.
Since $\cT$ consists of all the places not in $\cU$ or $\cS$, we thus
obtain from the product formula that
\begin{equation} 
 \begin{split}
   q^{d\deg P} \hhat(\beta) - C_4 \deg(Q)\log\deg(Q) - C_7 & \le
   \sum_{v \in \cT} - \log |\phi_P(\beta)|_v \\ & \le q^{d \deg P}
   \hhat(\beta) + C_8
\end{split}
\end{equation}
for all $P|Q$. Moreover, Lemma~\ref{inclusion-exclusion} gives
\begin{equation}\label{in-ex 2}
\begin{split}
  q^{d\deg P}\hhat(\beta)- C_4 \deg(Q)\log\deg(Q)-C_7 
  & \le \sum_{v \in \cT} \sum_{\Drin_P(y)=0} - \log |y-\beta|_{v} \\ &
  \leq q^{d \deg P} \hhat(\beta) + C_8.
\end{split}
\end{equation}

Using \eqref{in-ex 2}, we will bound
$$ \sum_{v \in \cT} \sum_{\substack{\Drin_P(y)=0 \\ \text {for
      $P | Q$ with $P \not = Q$}}} - \log |\beta - y|_{v}.$$
We compute the above sum via inclusion-exclusion (as we did in Lemma~\ref{inclusion-exclusion}) and we obtain
\begin{equation}
\label{Mobius}
\sum_{\substack{P|Q \\ P \not= Q}} -\mu\left(\frac{Q}{P}\right)
\sum_{v\in\cT} \sum_{\Drin_P(y)=0}-\log |\beta - y|_{v}.
    \end{equation}
    A simple computation using \eqref{in-ex 2} and \eqref{Mobius}
    shows that if $Q_1,\dots,Q_s$ are all the distinct irreducible
    factors of $Q$, then
\begin{equation}
\label{E:conclusion}
\begin{split}
 \sum_{v \in \cT} & \sum_{\substack{\phi_P(y) = 0 \\ \text { for
      $P | Q$ with $P \not = Q$}}} - \log |\beta - y|_{v} \\  
& \le \left(q^{d\deg(Q)} - q^{d\deg(Q)} \prod_{i=1}^s \left(
    1-q^{-d\deg(Q_i)} \right)\right)\hhat(\beta)\\
& \quad \quad  \quad + 2^s \left(C_4\deg(Q) \log\deg(Q) + C_7 + C_8 \right),
\end{split}
\end{equation}
due to the simple identity
$$\sum_{\substack{P|Q\\ P\not
    =Q}}-\mu\left(\frac{Q}{P}\right)q^{d\deg(P)}=q^{d\deg(Q)} -
q^{d\deg(Q)} \prod_{i=1}^s\left( 1-q^{-d\deg(Q_i)} \right)$$
and the
fact that there are $2^s-1$ nonzero terms in the outer sum from \eqref{Mobius}.  Using
\eqref{in-ex 2} for $P=Q$, we obtain
\begin{equation}
\label{E:conclusion 0}
\begin{split}
  \sum_{v \in \cT} \sum_{\phi_Q(y)=0} - \log |y-\beta|_{v}  \ge
  q^{d \deg Q} \hhat(\beta) - C_4 \deg(Q) \log\deg(Q) - C_7.
\end{split}
\end{equation}

But Lemma~\ref{one v} and equation \eqref{E:conclusion} imply that
\begin{equation}\label{Tq two}
\begin{split}
  \sum_{v \in \cT} & \sum_{\phi_Q(y)=0} - \log |y-\beta|_{v} \\ & \le \hhat(\beta) 
  \left( q^{d\deg(Q)} - q^{d\deg(Q)} \prod_{i=1}^s \left(
      1-q^{-d\deg(Q_i)} \right) \right) \\
  & \quad \quad +2^{s}\left(C_4\deg(Q) \log\deg(Q) + C_7 + C_8\right) + \deg Q.
\end{split}
\end{equation}

Since $\hhat(\beta)$ is positive, for $Q$ of large degree
we have
\begin{equation}
\label{E:contradiction 21}
\begin{split}
& q^{d \deg Q} \hhat(\beta) - C_4 \deg(Q) \log(\deg(Q)) - C_7 \\
&  > \hhat(\beta) \left( q^{d\deg(Q)} - q^{d\deg(Q)} \prod_{i=1}^s \left(
      1-q^{-d\deg(Q_i)} \right) \right) \\
& \quad \quad  + 2^{s}\left(C_4\deg(Q) \log\deg(Q) + C_7 +C_8\right) +  \deg Q
\end{split}
\end{equation}
because (recalling that $s$ is the number of distinct irreducible
factors of $Q$, which is smaller than $\deg Q$ as $\deg Q$ goes to infinity)
\begin{equation*}
\begin{split}
  q^{d \deg(Q)} & \prod_{i=1}^s \left(
      1-q^{-d\deg(Q_i)}\right) \hhat(\beta)
\ge \sqrt{q^{d\deg(Q)}}
      \hhat(\beta) \\
& >
  2^{s+1}\left(
    C_8 + \deg Q + C_4 \deg(Q) \log\deg(Q) + C_7\right).
\end{split}
\end{equation*}  
Inequalities \eqref{E:conclusion 0}, \eqref{Tq two} and \eqref{E:contradiction 21} give us a contradiction, which means that there exists
$v \in \cT$ such that $|\phi_Q(\beta)|_v < 1$ but
$|\phi_P(\beta)|_v = 1$ for all $P | Q$ with $P \not= Q$.  
\end{proof}

The following result is an immediate consequence of Theorem~\ref{precise order}.
\begin{cor}
\label{C:precise order}
  With the notation as in Theorem~\ref{precise order}, there exists a
  positive integer $N$ such that for all monic $Q\in\Fq[t]$ of degree at
  least $N$, there exists a place $v\in M_K$ of good reduction for
  $\phi$, such that $\beta$ is integral at $v$ and $\overline{\beta}$
  is a torsion point of order $Q$ for $\overline{\phi}$.
\end{cor}

\section{Further directions}\label{further}

Let $\psi: \bP^1 \lra \bP^1$ be a rational map of degree $d > 1$
defined over a number field $L$.  As with Drinfeld modules, one can
define canonical heights $\hhat_v$ for $\psi$, following Call and
Goldstine (\cite{CG}).  Pi{\~n}eiro, Szpiro, and Tucker (\cite{STP})
have proved that
$$
\hhat_v (\beta) = \int_{\bP^1(\bC_v)} \log |F|_v d \mu_{v, \psi},
$$
where $\hhat_v$ is the local canonical height for $\psi$ and
$\mu_{v,\psi}$ is an invariant measure associated to $\psi$ on
$\bP^1(\bC_v)$.  This generalizes an earlier formula due to Mahler
(\cite{mahler}). The invariant measure at an archimedean place $v$
was constructed by Lyubich (\cite{Lyubich}; see also \cite{Brolin} and
\cite{Mane, mane2}), while the invariant measure at nonarchimedean
places (which technically exists on the Berkovich space for
$\bP^1(\bC_v)$, as defined in \cite{berkovich}) was constructed by
Baker/Rumely (\cite{BR}), Chambert-Loir (\cite{CL}), and
Favre/Rivera-Letelier (\cite{FR1, FR2}).  These measures can also be
constructed at places of a function field over a finite field (see
\cite{Mat.Ann} for a treatment of these measures at places lying over
$v_{\infty}$). Thus, Theorem \ref{easiest} shows that a certain
integral of the function $\log|x - \beta|_v$ may be computed by
averaging this function over the torsion points of $\phi$.  Hence,
Theorem \ref{easiest} is a statement about equidistribution.  Note,
however, that statements about equidistribution are usually phrased in
terms of continuous functions (as is the case in \cite{Lyubich, BR,
  CL, FR1, FR2, Mat.Ann} for example). Theorem \ref{easiest} applies
to a function that has a logarithmic pole at an algebraic number.

Torsion points are merely inverse images of the point 0.  In light of
the results of \cite{ST} for rational functions over number fields,
one is led to make the following conjecture.

\begin{conjecture}\label{backwards}
  Let $\alpha$ be any point in $\bG_a(K)$.  Then, for any
  nontorsion $\beta \in \bK$ and any nonzero irreducible $F \in K[Z]$
  such that $F(\beta) = 0$, we have
\begin{equation*}
\begin{split}
  (\deg F) & \hhat(\beta) 
  = \sum_{v\in M_K} \lim_{\deg Q \to \infty}
  \frac{1}{{q^{d \deg Q}}} \sum_{\Drin_Q(y) = \alpha} \log | F(y) |_{v}.
\end{split}
\end{equation*}
\end{conjecture}

In \cite{BIR}, Baker, Ih, and Rumely show that local heights can also
be computed by averaging $\log |x - \beta|_v$ over {\it Galois orbits}
of torsion points of elliptic curves and of $\bG_m$ over a number field.
This is stronger, since the set of points $y$ such that $\psi^m(y) = e$
(where $e$ is the identity in the group) breaks into several Galois
orbits of torsion points.  The work of Baker, Ih, and Rumely leads us
to conjecture the following.

\begin{conjecture}\label{for-Ih}
  Let $(\alpha_n)_{n=1}^\infty$ be any nonrepeating sequence of
  torsion points in $\bG_a(\bK)$.  Then for any nontorsion $\beta \in
  \overline{K}$ and any nonzero irreducible $F \in K[Z]$ such that
  $F(\beta) = 0$, we have
  \begin{equation*}
  (\deg F) \hhat(\beta) 
  = \sum_{v\in M_K} \lim_{n \to \infty}
  \frac{1}{\#\left(\Gal(\bK/K)\cdot\alpha_n\right)}
  \sum\limits_{y
    \in{\Gal(\bK/K)\cdot \alpha_n}} \log |F(y)|_{v}.
\end{equation*}
\end{conjecture}

Arguing as in the proof of Theorem \ref{integral:1}, Conjecture
\ref{backwards} would yield the following corollary, which is
analogous to a theorem proved by Silverman (\cite{SilSiegel}) for
nonconstant morphisms of $\bP^1$ of degree greater than one over a
number field.

\begin{cor}\label{back-cor}
  Let $S$ be a finite set of places, let $\alpha$ be any point
  in $\bG_a(K)$, and let $\beta$ be any nontorsion point in
  $\bG_a(\bK)$. If Conjecture~\ref{backwards} holds, then there are finitely many $Q$ such that
  $\Drin_Q(\beta)$ is $S$-integral for $\alpha$.
\end{cor}

Similarly, Conjecture \ref{for-Ih} would imply the following
Corollary, which was proved by Baker, Rumely, and Ih (\cite{BIR}) for
elliptic curves and for $\bG_m$ over a number field.  Baker, Ih, and
Rumely have conjectured that an analog of this holds for arbitrary
nonconstant morphisms of degree greater than one of $\bP^1$ over a
number field.  

\begin{cor}
  Let $S$ be a finite set of places.  If Conjecture~\ref{for-Ih} holds, then for each nontorsion $\beta \in \bG_a(\bK)$ there are finitely many torsion points $\alpha$ that are
  $S$-integral for $\beta$.
\end{cor}

In the case of rational functions over number fields, the analog of
Conjecture \ref{backwards} has already been proved while the analog of
Conjecture \ref{for-Ih} is still just a conjecture (though some
special cases have been proved, as noted earlier).  In the case of
Drinfeld modules, on the other hand, a proof of Conjecture
\ref{for-Ih} does not seem far off.  Indeed, it would follow from a
combination of the Tate-Voloch conjecture (proved in \cite{TV}) with
suitable equidistribution results for continuous functions on the
Berkovich space for $\bP^1_K$.  Such results have been proved in the
number field case by Baker/Rumely (\cite{BR}), Chambert-Loir
(\cite{CL}), and Favre/Rivera-Letelier (\cite{FR1, FR2}).  M.~Baker
informs us that the techniques used in these proofs should also work
over function fields over finite fields.

It is less clear how one might proceed towards a proof of Conjecture
\ref{backwards}.  The arguments involving linear forms in logarithms
at the infinite places work exactly the same, but the proof of Lemma
\ref{L:nice trick} does not seem to work when 0 is replaced by an
arbitrary point on $\bG_a(\bK)$.  Nevertheless, we are tempted to ask
whether something even more general than Corollary~\ref{back-cor} is
true.

\begin{conjecture}
\label{Siegel-Drinfeld}
  Let $S$ be a finite set of places, let $\alpha$ be any point
  in $\bG_a(\overline{K})$, and let $\Gamma$ be any finitely generated
  $\phi$-submodule of $\bG_a(\overline{K})$. Then there are then finitely many $\gamma \in \Gamma$ such that
  $\gamma$ is $S$-integral for $\alpha$.
\end{conjecture}

This would say the same thing for finitely generated
$\Drin$-submodules of $\bG_a(\overline{K})$ what Siegel's theorem says
for finitely generated $\bZ$-submodules of $\bG_m(\overline{K})$.
Theorem \ref{integral:1} shows that our Conjecture~\ref{Siegel-Drinfeld} holds when $\Gamma$ is a cyclic $\Drin$-submodule and $\alpha = 0$
(note that is also a special case of Corollary~\ref{back-cor}).  We can
easily generalize this to a result that applies for any
$\alpha\in\phi_{\tor}(K)$ (when $\Gamma$ is a cyclic $\phi$-submodule),
as we shall see next. We also mention that in \cite{our_siegel} we proved Conjecture~\ref{Siegel-Drinfeld} under the assumption that $K$ has only one infinite place, and in addition assuming that Bosser's result on linear forms in logarithms holds also at finite places. While the second assumption is believed to be true by experts in the field, we also believe that the assumption from our paper \cite{our_siegel} of having only one infinite place in $K$ could be removed, and so, Conjecture~\ref{Siegel-Drinfeld} would hold in the generality we stated it here.

Following the same reasoning we used in deducing
Theorem~\ref{integral:1} from Theorem~\ref{easiest} (which in turn was
deduced from Corollary~\ref{alternative local height computation}), we
can prove Conjecture~\ref{Siegel-Drinfeld} in the case
$\alpha\in\phi_{\tor}(K)$ and $\Gamma$ is a cyclic $\phi$-submodule
generated by a nontorsion point $\beta$, from the following
Proposition.

\begin{proposition}
\label{torsion-Siegel}
If $\beta$ is a non-torsion point, and $\alpha\in K$ is a torsion point,
then for every place $w\in M_{K(\beta)}$, we have
$$\hhat_w(\beta)=\lim_{\deg Q\to\infty}\frac{\log
  |\phi_Q(\beta)-\alpha |_w}{q^{d\deg Q}}.$$
\end{proposition}

We noted earlier that Corollary~\ref{back-cor} would follow easily
from Conjecture~\ref{backwards} via the arguments used in the
derivation of Theorem~\ref{integral:1} from Theorem~\ref{easiest}.
Now, reasoning as in the proof of Theorem~\ref{easiest}, we see that
Proposition~\ref{torsion-Siegel} yields Conjecture~\ref{backwards} in
the case where $\alpha$ is torsion because (after reducing to the case
that $F(Z)=\prod_{i=1}^n (Z-\theta_i)$ is separable, as we did in the proof of
Theorem~\ref{easiest}) we have
$$\sum_{i=1}^n\sum_{\phi_Q(y)=\alpha} \log
|y-\theta_i|_v=n\cdot\left(\sum_{\substack{w|v\\w\in
      M_{K(\beta)}}}\log|\phi_Q(\beta)-\alpha|_w\right)- n
\log |\gamma_Q|_v$$
for each $v\in M_K$ (where $\gamma_Q$ is as before
the leading coefficient of $\phi_Q$). The above equality follows
readily from Proposition~\ref{torsion-Siegel} using the arguments that
appeared in the proof of Theorem~\ref{easiest}, after noting the
bijection between the two sets
$$\{y\text{ : }\phi_Q(y)=\alpha\}\text{ and }\{\alpha_0+y\text{ :
}\phi_Q(y)=0\},$$
where $\alpha_0$ is a fixed solution for
$\phi_Q(x)=\alpha$. Therefore all we need to do is prove
Proposition~\ref{torsion-Siegel}.  We provide a sketch of a proof
below.

\begin{proof}[Proof of Proposition~\ref{torsion-Siegel}.]
  If $\beta$ is not in the filled Julia set for $\phi$ at $w$, then
  $$|\phi_Q(\beta)-\alpha|_w=|\phi_Q(\beta)|_w$$
  for polynomials $Q$
  of large degree. Thus, in this case, the conclusion of
  Proposition~\ref{torsion-Siegel} is immediate.
  
  Assume from now on that $\beta$ is in the filled Julia set for
  $\phi$ at $w$. Then we need to show that
\begin{equation}
\label{t-S-1}
\lim_{\deg Q\to\infty}\frac{\log |\phi_Q(\beta)-\alpha|_w}{q^{d\deg Q}} = 0.
\end{equation}
Let $\epsilon>0$. Since $|\phi_Q(\beta)|_w$ is bounded as $\deg
Q\to\infty$, it follows that
\begin{equation}
\label{t-S-2}
\limsup_{\deg Q\to\infty}\frac{\log |\phi_Q(\beta)-\alpha|_w}{q^{d\deg Q}} < \epsilon.
\end{equation}
It will suffice to show that also
\begin{equation}
\label{t-S-3}
\liminf_{\deg Q\to\infty}\frac{\log |\phi_Q(\beta)-\alpha|_w}{q^{d\deg Q}} > -\epsilon.
\end{equation}
If $w$ lies over $v_{\infty}$, then an argument almost identical to
the one used in Proposition~\ref{P:first iteration} will give rise to
\eqref{t-S-3}.  As with Proposition~\ref{P:first iteration}, the key
ingredient is an application of the lower bound for linear forms in
logarithms as provided by Proposition~\ref{from Bosser}.

If $w$ does not lie over $v_{\infty}$, then a slight modification of
our argument from Proposition~\ref{P:first finite place iteration} may
be used to prove \eqref{t-S-3}. For this we use the following claim.

\begin{claim}
\label{nice trick again}
Let $w$ be a place that does not lie over $v_\infty$.  Then there
exists a positive constant $C_w$ such that for every $Q\in A$, we have
$|\phi_Q(\beta)-\alpha|_w\ge C_w|Q|_w$.
\end{claim}

\begin{proof}[Proof of Claim~\ref{nice trick again}.]
  Let $R\in A$ be a non-constant polynomial such that
  $\phi_R(\alpha)=0$ (we recall that $\alpha\in\phi_{\tor}$). If for
  every $Q\in A$, we have $|\phi_Q(\beta)-\alpha|_w\ge |R|_w$, then
  the conclusion of Claim~\ref{nice trick again} is immediate.
  Therefore, assume there exists $Q\in A$ such that
  $|\phi_Q(\beta)-\alpha|_w<|R|_w$. For each such $Q$,
  Lemma~\ref{small-x} yields $
  |\phi_R(\phi_Q(\beta)-\alpha)|_w=|R|_w\cdot
  |\phi_Q(\beta)-\alpha|_w,$ which implies that
\begin{equation}
\label{t-S-4}
|\phi_{RQ}(\beta)|_w = |R|_w\cdot |\phi_Q(\beta)-\alpha|_w,
\end{equation}
since $\phi_R(\alpha) = 0$.  By Lemma~\ref{L:nice trick}, there exists
$C_w>0$ such that for every polynomial $Q'$, we have
$|\phi_{Q'}(\beta)|_w\ge C_w |Q'|_w$. Hence, in particular, we have
$|\phi_{RQ}(\beta)|_w\ge C_w\cdot |RQ|_w$.  Applying \eqref{t-S-4}, we
obtain $|\phi_Q(\beta)-\alpha|_w\ge C_w|Q|_w$, as desired.
\end{proof}
Since $\lim_{\deg Q\to\infty}\frac{\log |Q|_w}{q^{d\deg Q}}=0$,
Proposition~\ref{torsion-Siegel} now follows immediately from
Claim~\ref{nice trick again}.
\end{proof}

\def\cprime{$'$} \def\cprime{$'$} \def\cprime{$'$} \def\cprime{$'$}
\providecommand{\bysame}{\leavevmode\hbox to3em{\hrulefill}\thinspace}
\providecommand{\MR}{\relax\ifhmode\unskip\space\fi MR }
\providecommand{\MRhref}[2]{%
  \href{http://www.ams.org/mathscinet-getitem?mr=#1}{#2}
}
\providecommand{\href}[2]{#2}


\end{document}